\documentstyle[12pt]{article}

\def\bl{\begin{eqnarray}}
\def\el{\end{eqnarray}}
\def\bll{\begin{eqnarray*}}
\def\ell{\end{eqnarray*}}
\def\su{\sum\limits}

\begin{document}

\title{A Generalization of Monotonicity Condition and Applications\thanks{%
Supported in part by Natural Science Foundation of China under grant number
10471130.}}
\author{D. S. Yu and S. P. Zhou}
\date{}
\maketitle

\pagenumbering{arabic}

\begin{quote}
{\small \ {\bf ABSTRACT: {\rm In the present paper, we introduce a new class
of sequences called  $NBVS$ to generalize $GBVS$, essentially extending
monotonicity from \textquotedblleft one sided" to \textquotedblleft two
sided", while some important classical results keep true. }}}
\end{quote}

\begin{center}
{\small 2000 Mathematics Subject Classification: 42A20, 42A32.}
\end{center}

\begin{flushleft}
{\bf \S 1. Introduction}
\end{flushleft}

{\normalsize It is well known that there are a great number of interesting
results in Fourier analysis established by assuming monotonicity of
coefficients, and many of them have been generalized by loosing the
condition to quasi-monotonicity, $O$-regularly varying quasi-monotonicity,
etc.. Generally speaking, it has become an important topic how to generalize
monotonicity. }

{\normalsize Recently, Leindler [5] defined a new class of sequences named
as {\it sequences of rest bounded variation} and briefly denoted by $RBVS$,
keeping some good properties of decreasing sequences. A null sequences ${\bf %
C} :=\{c_n\}$ is of rest bounded variation, or ${\bf C}\in RBVS$, if $%
c_n\rightarrow 0$ and for any $m\in N$ it holds that
\[
\sum\limits_{n=m}^\infty|\Delta c_n|\leq K({\bf C})c_m,
\]
where $\Delta c_n=c_n-c_{n+1},$ and $K({\bf C})$ denotes a constant only
depending on ${\bf C}$. By the definition of $RBVS$, it clearly yields that
if ${\bf C}\in RBVS$, then for any $n\geq m$, it holds that
\[
c_n\leq K({\bf C})c_m.
\]
Denote by $MS$ the monotone decreasing sequence and $CQMS$ the classic
quasi-monotone decreasing sequences\footnote{${\bf C}=\{c_{n}\}\in CQMS$
means that there is an $\alpha\geq 0$ such that $c_{n}/n^{\alpha}$ is
decreasing.}, then it is obvious that
\begin{eqnarray*}
MS\subset RBVS\cap CQMS.
\end{eqnarray*}
}

{\normalsize Leindler [6] proved that the class $CQMS$ and $RBVS$ are not
comparable. Very recently, Le and Zhou [3] suggested the following  new
class of sequences to include both $RBVS$ and $CQMS$: }

{\normalsize {\bf Definition 1.}\quad {\it Let ${\bf C} :=\{c_n\}$ be a
sequence satisfying $c_n\in M(\theta_0):=\big
 \{z:|\arg z|\leq \theta_0\big\}$ for some $\theta_0\in [0,\pi/2)$ and $%
n=1,2,\cdots.$  If there is an integer $N_{0}\geq 1$ such that
\[
\sum\limits_{k=m}^{2m}|\Delta c_k|\leq K({\bf C})\max_{m\leq n<m+N_{0}}|c_n|
\]
holds for all $m=1,2,\cdots,$ where here and throughout the paper, $K({\bf C}%
)$ always indicates a positive constant only depending on ${\bf C}$, and its
value may be different even in the same line, then we say ${\bf C}\in GBVS.$}
}

{\normalsize Le and Zhou [3] proved that
\begin{eqnarray*}
RBVS\cup CQMS\subset GBVS.
\end{eqnarray*}
}

{\normalsize If $\{c_n\}$ is a nonnegative sequence tending to zero for a
special case $N_{0}=1$, then $GBVS$ can be restated as }

{\normalsize {\bf Definition 1$^{\prime }$.} \quad {\it Let ${\bf C}
:=\{c_n\}$ be a nonnegative sequence tending to zero, if
\[
\sum\limits_{n=m}^{2m}|\Delta c_n|\leq K({\bf C})c_m
\]
holds for all $m=1,2,\cdots,$ then we say ${\bf C}\in GBVS.$} }

{\normalsize For convenience, it will refer to $GBVS$ in Definition $%
1^{\prime }$ (or $GBVS$ in Definition 1 for $N_{0}=1$) when we mention $GBVS$
in the present paper later. }

{\normalsize After some interesting investigations, Zhou and Le [12] comes
to reach a conclusion by saying that ``{\it In any sense, monotonicity and
`rest bounbded variation' condition are all `one sided' monotonicity
condition, that is,  a positive sequence ${\bf b}=\{b_n\}$ under any of
these conditions satisfies $b_n\leq Cb_k$  for $n\geq k: b_n$ can be
controlled by one factor $b_k$. But for $\{b_n\}\in GBVS$, one can
calculate, for $k\leq n\leq 2k,$
\[
b_n=\sum\limits_{j=n}^{2k}\Delta
b_j+b_{2k+1}\leq\sum\limits_{j=n}^{2k}|\Delta b_j|+b_{2k+1}  \leq K({\bf b}%
)b_k+b_{2k+1},
\]
and this can be actually regarded as a `two sided' monotonicity: $b_n$ is
controlled not  only by $b_k$ but also by $b_{2k+1}$. Therefore, the
essential point of $GBV$ condition is to  extend monotonicity from `one
sided' to  `two sided'.}" However, by the definition of $GBVS$, we finally
realize to  deduce that
\begin{eqnarray*}
b_{2k+1}\leq \sum\limits_{j=k}^{2k}|\Delta b_j|+|b_{k}|\leq K({\bf b})b_{k},
\end{eqnarray*}
hence, $b_n\leq K({\bf b})b_k$ for all $k\leq n\leq 2k$ means that $GBV$
condition is still a ``one sided" condition in some local sense. Therefore
the object of the present paper is to suggest a new  class of sequence named
as $NBVS$ to include $GBVS$ (for the special case $N_{0}=1$)  essentially
extending monotonicity  from ``one sided" to ``two sided", while some
important classical results keep true. }

\begin{flushleft}
{\normalsize {\bf \S 2. Main Results} }
\end{flushleft}

{\normalsize {\bf Definition 2.}\quad {\it Let ${\bf C} :=\{c_n\}$ be a
sequence satisfying $c_n\in M(\theta_0):=\big
 \{z:|\arg z|\leq \theta_0\big\}$ for some $\theta_0\in [0,\pi/2)$ and $%
n=1,2,\cdots.$  If
\begin{eqnarray*}
\sum\limits_{n=m}^{2m}|\Delta c_n|\leq K({\bf C})\big(|c_m|+|c_{2m}|\big)
\end{eqnarray*}
holds for all $m=1,2,\cdots,$ then we say ${\bf C}\in NBVS.$} }

{\normalsize The following result implies that $NBVS$ is an essential
generalization of $GBVS$. }

{\normalsize {\bf Theorem 1.}\quad {\it Let ${\bf C}:=\{c_n\}$, if ${\bf C}%
\in GBVS$ $($for the special case $N_{0}=1$$)$, then ${\bf C}\in NBVS$, but
the converse is not true, that is, there exists a null sequence ${\bf C}$
such that ${\bf C}\in NBVS$, but ${\bf C}\not\in GBVS$.} }

{\normalsize Next we will prove that some important classical results still
keep true by applying this new condition. }

{\normalsize Let $C_{2\pi}$ be the space of all complex valued continuous
functions $f(x)$ of period $2\pi$ with the norm
\[
\|f\|:=\max\limits_{-\infty<x<\infty}|f(x)|.
\]
Given a trigonometric series $\sum\limits_{k=-\infty}^{\infty}c_ke^{ikx}:=%
\lim\limits_{n\rightarrow\infty} \sum\limits_{k=-n}^{n}c_ke^{ikx},$ write
\[
f(x)=\sum\limits_{k=-\infty}^{\infty}c_ke^{ikx}
\]
at any point $x$ where the series converges, and denote its $n$th partial
sum $S_n(f,x)$ by $\sum\limits_{k=-n}^{n}c_ke^{ikx}.$ }

{\normalsize {\bf Theorem 2.}\quad {\it Let ${\bf C}:=\{c_n\}_{n=0}^\infty$
be a complex sequence satisfying
\begin{eqnarray*}
c_n+c_{-n}\in M(\theta_0),\;n=1,2,\cdots
\end{eqnarray*}
for some $\theta_0\in [0,\pi/2)$, and ${\bf C}\in NBVS$, then the necessary
and sufficient conditions for $f\in C_{2\pi}$ and $\lim\limits_{n\rightarrow%
\infty}\|f-S_n(f)\|=0$ are that
\begin{eqnarray*}
\lim\limits_{n\rightarrow\infty}nc_n=0 \hspace{.4in}(1)
\end{eqnarray*}
and
\begin{eqnarray*}
\sum\limits_{n=1}^\infty|c_n+c_{-n}|<\infty.\hspace{.4in}(2)
\end{eqnarray*}
} }

{\normalsize As an application of Theorem 2, we have }

{\normalsize {\bf Theorem 3.}\quad {\it Let ${\bf b}=\{b_n\}_{n=1}^\infty$
be a nonnegative sequences. If ${\bf b}\in NBVS$, then a the necessary and
sufficient conditions either for the uniform convergence of series $%
\sum\limits_{k=1}^\infty b_n\sin nx$, or the continuity of its sum function $%
f(x)$, is that $\lim\limits_{n\rightarrow\infty}nb_n=0.$ } }

{\normalsize Theorem 3 can be derived by applying Theorem 2 by using exactly
the same technique as [3] did, we omit the details. We remark that Theorem 3
is first established for decreasing real sequences $\{b_n\}$ by Chaundy and
Jolliffe [2], for $\{b_n\}\in CQMS$ by Nurcombe [7], for $\{b_n\}\in RBVS$
by Leindler [5], and for general $\{c_n\}\in GBVS$ by Le and Zhou [3]. }

{\normalsize Denote by $E_n(f)$ the best approximation of $f$ by
trigonometric polynomials of degree $n$. Then }

{\normalsize {\bf Theorem 4}\quad {\it Let $\left\{c_{n}\right\}_{n=0}^%
\infty\in NBVS$ and $\left\{c_{n}+c_{-n}\right\}_{n=0}^\infty\in NBVS$, and
\[
f(x)=\sum\limits_{n=-\infty}^\infty c_{n}e^{inx}.
\]
Then $f\in C_{2\pi}$ if and only if
\[
\lim\limits_{n\rightarrow\infty}n c_{n}=0.\hspace{.4in}(1^{\prime })
\]
and
\[
\sum\limits_{n=1}^\infty|c_{n}+c_{-n}|<\infty. \hspace{.4in}(2^{\prime })
\]
Furthermore, if $f\in C_{2\pi}$, then
\[
E_n(f)\sim\max\limits_{1\leq k\leq n}k\left(|c_{n+k}|+|c_{-n-k}|\right)
+\max\limits_{k\geq 2n+1}k|c_{k}-c_{-k}|+\sum\limits_{k=2n+1}^{\infty}
|c_{k}+c_{-k}|.
\]
} }

{\normalsize Theorem 4 was established by Zhou and Le [11] for $GBVS$,
readers can find further related information in [11]. }

{\normalsize Let
\[
f(x) = \sum_{n=0}^{\infty}c_{n}\cos nx.
\]
{\rm The following theorem is an interesting application to a hard problem
in classical Fourier analysis. }}

{\normalsize {\rm {\bf Corollary 1.}\quad {\it Let $\{c_{n}\}\in\mbox{\rm
NBV}$ be a real sequence. If $f\in C_{2\pi}$ and
\[
\sum_{k=n+1}^{2n}c_{k}=O\left(\max_{1\leq k\leq n}kc_{n+k}\right),
\]
then
\[
\|f-S_{n}(f)\| = O(E_{n}(f)).
\]%
} }}

{\normalsize {\rm The original form of Corollary 1 for $GBVS$ could be found
in [11]. }}

{\normalsize {\rm Let $L_{2\pi}$ be the space of all complex valued
integrable functions $f(x)$ of period $2\pi$ with the norm
\[
\|f\|_L=\int_{-\pi}^\pi|f(x)|dx.
\]
Denote the Fourier series of $f\in L_{2\pi}$ by $\sum\limits_{k=-\infty}^{%
\infty}\hat{f}(k)e^{ikx}$. }}

{\normalsize {\rm {\bf Theorem 5.}\quad {\it Let $f(x)\in L_{2\pi}$ be a
complex valued function. If the Fourier coefficients of $f$ satisfying that
both $\{\hat{f}(n)\}_{n=0}^{+\infty}\in NBVS$ and $\{\hat{f}%
(-n)\}_{n=0}^{+\infty}\in NBVS$, then
\[
\lim\limits_{n\rightarrow\infty}\|f-S_n(f)\|_L=0\;\mbox{if and only if}\;
\lim\limits_{n\rightarrow\infty}\hat{f}(n)\log |n|=0.
\]%
} As a special case, we have }}

{\normalsize {\rm {\bf Corollary 2.} \quad {\it Let $f(x)\in L_{2\pi}$ be a
real valued function. If $\{\hat{f}(n)\}_{n=0}^{+\infty}\in NBVS$, then
\[
\lim\limits_{n\rightarrow\infty}\|f-S_n(f)\|_L=0\;\mbox{if and only if}\;
\lim\limits_{n\rightarrow\infty}\hat{f}(n)\log n=0.
\]%
} }}

{\normalsize {\rm Let $E_n(f)_L$ be the best approximation of a complex
valued function $f\in L_{2\pi}$ by trigonometric polynomials of degree $n$
in the integral metric, that is,
\[
E_n(f)_L:=\inf_{c_k}\left\|f-\sum\limits_{k=-n}^{n}c_ke^{ikx}\right\|_L.
\]
We establish the following $L^1-$approximation theorem. }}

{\normalsize {\rm {\bf Theorem 6.}\quad {\it Let $f(x)\in L_{2\pi}$ be a
complex valued function, $\{\psi_n\}$ be a decreasing sequence tending to
zero with that $\psi_n=O(\psi_{2n})$. If both $\{\hat{f}(n)\}_{n=0}^{+%
\infty}\in NBVS$ and $\{\hat{f}(-n)\}_{n=0}^{+\infty}\in NBVS$, then
\[
\lim\limits_{n\rightarrow\infty}\|f-S_n(f)\|_L=O(\psi_{n})
\]
if and only if
\[
E_n(f)_L=O(\psi_{n})\;\mbox{and}\;\hat{f}(n)\log |n|=O(\psi_{|n|}).
\]%
} }}

{\normalsize {\rm Theorem 5 and Theorem 6 were established for $O$-regularly
varying quasi-monotone sequences by Xie and Zhou [10], and for $GBVS$ by Le
and Zhou [4]. }}

\begin{flushleft}
{\normalsize {\rm {\bf \S 3. Proof of Results} }}

{\normalsize {\rm {\bf \S 3.1. Proof of Theorem 1} }}
\end{flushleft}

{\normalsize {\rm Let ${\bf C}:=\{c_n\}$. By the definition of $GBVS$ and $%
NBVS$, it is obvious that ${\bf C}\in GBVS$ implies ${\bf C}\in NBVS.$ }}

{\normalsize {\rm Set $k_j=2^j,j=1,2,\cdots.$ Define
\begin{eqnarray*}
c_n:=\left\{%
\begin{array}{ll}
\frac{1}{k_j^2}, & j\;\mbox{is even}, \\
\frac{1}{k_j^3}, & j\;\mbox{is odd}%
\end{array}%
\right.
\end{eqnarray*}
for $k_j\leq n<k_{j+1}$ and $c_1=1$. Let $2^j\leq m<2^{j+1}$ (without loss
of generality, assume that $j>2$). Then
\begin{eqnarray*}
\sum\limits_{k=m}^{2m}|\Delta
c_k|&=&|c_{2^{j+1}-1}-c_{2^{j+1}}|=|c_{2^j}-c_{2^{j+1}}| \\
&\leq &\left\{%
\begin{array}{ll}
\frac{1}{2^{2j+2}}=c_{2m}, & j \;\mbox{is odd,} \\
\frac{1}{2^{2j}}=c_{m}, & j \;\mbox{is even.}%
\end{array}%
\right.
\end{eqnarray*}
Therefore, ${\bf C}\in NBVS.$ On the other hand, Let $m=2^{2n+1}$, we have
\[
\sum\limits_{k=m}^{2m}|\Delta c_k|=\frac{1}{(2m)^2}-\frac{1}{m^3}\geq\frac{1%
}{8m^2},
\]
that means,
\[
\frac{\sum\limits_{k=m}^{2m}|\Delta c_k|}{c_m}\geq \frac{m}{8}.
\]
In other words, there is no constant $K({\bf C})$ such that
\[
\sum\limits_{k=m}^{2m}|\Delta c_k|\leq K({\bf C})c_m
\]
holds for all $m\geq 1.$ Therefore, ${\bf C}\not\in GBVS.$ }}

\begin{flushleft}
{\normalsize {\rm {\bf \S 3.2. Proof of Theorem 2} }}
\end{flushleft}

{\normalsize {\rm {\bf Lemma 1 (Xie and Zhou [9]).}\quad {\it Let ${\bf C}%
:=\{c_n\}_{n=0}^\infty$ be a complex sequence satisfying
\begin{eqnarray*}
c_n+c_{-n}\in M(\theta_0),\;n=1,2,\cdots
\end{eqnarray*}
for some $\theta_0\in [0,\pi/2)$. Then $f\in C_{2\pi}$ implies that
\[
\sum\limits_{n=1}^\infty|c_n+c_{-n}|<\infty.
\]%
} }}

{\normalsize {\rm {\bf Lemma 2.}\quad {\it Let ${\bf C}:=\{c_n\}_{n=0}^\infty
$ satisfy all conditions of Theorem $2$. Then $\lim\limits_{n\rightarrow%
\infty}\|f-S_n(f)\|=0$ implies that
\[
\lim\limits_{n\rightarrow\infty}nc_n=0.
\]%
} }}

{\normalsize {\rm {\bf Proof.} A standard calculation yields that
\[
S_{4n}(f,x)-S_n(f,x)=\sum\limits_{k=n+1}^{4n}c_k\left(e^{ikx}-e^{-ikx}%
\right)+ \sum\limits_{k=n+1}^{4n}(c_k+c_{-k})e^{-ikx}.
\]
Let $x_0=\frac{\pi}{8n}$, then
\begin{eqnarray*}
\frac{1}{2} \sum\limits_{k=n+1}^{4n} \mbox{Re} c_k&\leq&2\sum%
\limits_{k=n+1}^{4n}\mbox{Re} c_k\sin kx_0=\left|\sum\limits_{k=n+1}^{4n}%
\mbox{Re} c_k\left(e^{ikx_0}-e^{-ikx_0}\right)\right|  \nonumber \\
&\leq&\left|\sum\limits_{k=n+1}^{4n}
c_k\left(e^{ikx_0}-e^{-ikx_0}\right)\right|  \nonumber \\
&\leq& \|S_{4n}(f)-S_n(f)\|+\sum\limits_{k=n+1}^{4n}|c_k+c_{-k}|.  \label{C2}
\end{eqnarray*}
If $a,b\in K(\theta_0)$ for some $\theta_0\in [0,\pi/2)$, then $a+b\in
K(\theta_0)$. Hence the condition $\lim\limits_{n\to\infty}\|f-S_{n}(f)\|=0$
(consequently $c_{n}\to 0$ as $n\to\infty$) implies $|c_n|\leq M(\theta_0)%
\mbox{Re}c_n$ for all $n\geq 1$. By the definition of $NBVS$, for $n\leq
k\leq 2n$, it holds that
\begin{eqnarray*}
|c_{2n}|\leq \sum\limits_{j=k}^{2n-1}|\Delta c_j|+|c_k|\leq
\sum\limits_{j=k}^{2k}|\Delta c_j|+|c_k|\leq K({\bf C})\left(|c_k|+|c_{2k}|%
\right).\hspace{.4in}(3)
\end{eqnarray*}
Thus
\begin{eqnarray*}
n|c_{2n}|&\leq& K({\bf C})\sum\limits_{k=n+1}^{2n}\left(|c_k|+|c_{2k}|%
\right)\leq K({\bf C})\sum\limits_{k=n+1}^{4n}|c_k| \\
&\leq&K({\bf C},\theta_0)\sum\limits_{k=n+1}^{4n}\mbox{Re} c_k.
\end{eqnarray*}
By Lemma 1 with all the above estimates, we get $\lim\limits_{n\rightarrow%
\infty}nc_{2n}=0$. A similar discussion also yields $\lim\limits_{n%
\rightarrow\infty}nc_{2n+1}=0$. Lemma 2 is proved. }}

{\normalsize {\rm {\bf Proof of Theorem 2.} We will follow the line of proof
in [3], only some necessary modifications will be mentioned. }}

{\normalsize {\rm {\bf Necessity.} Applying Lemma 1 and Lemma 2, we
immediately have (1) and (2). }}

{\normalsize {\rm {\bf Sufficiency.} It is not difficult to see that, under
the conditions (1) and (2), $\{S_n(f,x)\}$ is a Cauchy sequence for each $x$%
, consequently it converges at each $x$. Now we need only to show that
\begin{eqnarray*}
\lim\limits_{n\rightarrow\infty}\left\|\sum\limits_{k=n}^\infty
\left(c_ke^{ikx}+c_{-k}e^{-ikx}\right)\right\|=0  \label{C4}
\end{eqnarray*}
in this case. In view of (1) and (2), for any given $\varepsilon>0$, there
is a $n_0$ such that for all $n\geq n_0$, it holds that
\begin{eqnarray*}
\sum\limits_{k=n}^\infty|c_k+c_{-k}|<\varepsilon \hspace{.4in}(4)
\end{eqnarray*}
and
\begin{eqnarray*}
n|c_n|<\varepsilon.\hspace{.4in}(5)
\end{eqnarray*}
Let $n\geq n_0$, set
\begin{eqnarray*}
\sum\limits_{k=n}^\infty \left(c_ke^{ikx}+c_{-k}e^{-ikx}\right)=
\sum\limits_{k=n}^\infty(c_k+c_{-k})e^{-ikx}+2i\sum\limits_{k=n}^\infty
c_k\sin kx:=I_1(x)+2iI_2(x).
\end{eqnarray*}
}}

{\normalsize {\rm By (4), we get
\[
I_1(x)<\varepsilon.
\]
For $x=0$ and $x=\pi$, we have $I_2(x)=0$. Therefore, we may restrict $x\in
(0,\pi).$ Take $N:=[1/x]$ and set\footnote{%
When $N\leq n$, the same argument as in estimating $J_2$ can be applied to
deal with $I_2(x)=\sum\limits_{k=n}^\infty c_k\sin kx$ directly.}
\[
I_2(x)=\sum\limits_{k=n}^{N-1}c_k\sin kx+\sum\limits_{k=N}^\infty c_k\sin
kx:= J_1(x)+J_2(x).
\]
From (5) with $N=[1/x]$, it follows that
\[
|J_1(x)|\leq \sum\limits_{k=n}^{N-1}|c_k\sin kx|\leq
x\sum\limits_{k=n}^{N-1}k|c_k| <x(N-1)\varepsilon\leq \varepsilon.
\]
Set $D_n(x)=\sum\limits_{k=1}^n\sin kx,$ it is well known that $|D_n(x)|\leq%
\frac{\pi}{x}.$ By Abel's transformation, we have
\begin{eqnarray*}
|J_2(x)|&=&\sum\limits_{k=N}^\infty|\Delta c_k||D_k(x)|+c_N|D_{N-1}(x)|\leq
Mx^{-1}\left(\sum\limits_{k=N}^\infty|\Delta c_k|+|c_N|\right)  \nonumber \\
&\leq&K({\bf C})\left(x^{-1}\sum\limits_{k=N}^\infty|\Delta
c_k|+N|C_N|\right)\leq K({\bf C})\left(x^{-1}\sum\limits_{k=N}^\infty|\Delta
c_k|+\varepsilon\right). \hspace{.4in}(6)
\end{eqnarray*}
Since ${\bf C}\in NBVS$, we calculate that
\begin{eqnarray*}
\sum\limits_{k=N}^\infty|\Delta c_k| &\leq&\sum\limits_{k\geq \log N/\log
2}\sum\limits_{j=2^k}^{2^{k+1}-1} |\Delta c_k|+\sum\limits_{k=N}^{2N}|\Delta
c_k|  \nonumber \\
&\leq&K({\bf C})\left(\sum\limits_{ k\geq \log N/\log
2}(|c_{2^k}|+|c_{2^{k+1}}|)+|c_N|+|c_{2N}|\right)  \nonumber \\
&\leq&K({\bf C})\left(\sum\limits_{ k\geq \log N/\log 2} 2^{-k}
\left(2^k|c_{2^k}|+2^{k+1}|c_{2^{k+1}}|\right)+
\left(N|c_N|+2N|c_{2N}|\right)N^{-1}\right)  \nonumber \\
&\leq&K({\bf C})\varepsilon\left(\sum\limits_{ k\geq \log N/\log
2}2^{-k}+N^{-1}\right)\leq K({\bf C})\varepsilon N^{-1}. \hspace{.4in}(7)
\end{eqnarray*}
Combining (6) and (7) yields that $|J_2(x)|\leq  K({\bf C})\varepsilon,$ and
altogether it follows that
\[
\lim\limits_{n\rightarrow\infty}\|f-S_n(f)\|=0.
\]
}}

\begin{flushleft}
{\normalsize {\rm {\bf \S 3.3. Proof of Theorem 4} }}
\end{flushleft}

{\normalsize {\rm {\bf Lemma 3.}\quad {\it Let $\left\{c_{n}\right\}_{n=0}^%
\infty\in NBVS$, $\left\{c_{n}+c_{-n}\right\}_{n=0}^\infty\in NBVS$, and
\[
f(x)=\sum\limits_{n=-\infty}^\infty c_{n}e^{inx}.
\]
Suppose $f\in C_{2\pi}$, then
\[
\max\limits_{k\geq 1}k|c_{\pm(n+k)}|=O(E_n(f)).
\]
} }}

{\normalsize {\rm {\bf Proof.} Let $t_n^*(x)$ be the trigonometric
polynomials of best approximation of  degree $n$, then from the obvious
equality
\[
\frac{1}{2\pi}\int_{-\pi}^\pi\left|e^{\pm
i(n+1)x}\left(\sum\limits_{k=0}^{N-1}e^{\pm ikx}  \right)^2\right|dx=N
\]
we get
\[
\sum\limits_{k=1}^N\left(kc_{\mp (n+k)}+(N-k)c_{\mp (n+N+k)}\right)
\]
\[
=\left|\frac{1}{2\pi}\int_{-\pi}^\pi  (f(x)-t_n^*(x))e^{\pm
i(n+1)x}\left(\sum\limits_{k=0}^{N-1}e^{\pm ikx}  \right)^2dx\right|\leq
NE_n(f).
\]
Thus
\begin{eqnarray*}
NE_n(f)&\geq&\mbox{Re}\left(\sum\limits_{k=1}^N
\left(kc_{\mp(n+k)}+(N-k)c_{\mp (n+N+k)}\right)\right)  \nonumber \\
&\geq&\sum\limits_{k=1}^N\left(k\mbox{Re}c_{n+k}+(N-k)\mbox{Re}%
c_{n+N+k}\right)  \nonumber \\
&\geq&\sum\limits_{k=1}^Nk\mbox{Re}c_{n+k}.\hspace{.4in}(8)
\end{eqnarray*}
For any $\frac{N}{2}\leq j\leq N$, it follows from the definition of $NBVS$
that
\[
|c_{n+N}|\leq\sum\limits_{k=j}^{N-1}|\Delta
c_{n+k}|+|c_{n+j}|\leq\sum\limits_{k=j}^{2j}|\Delta c_{n+k}|+|c_{n+j}|
\]
\[
\leq K \left(|c_{n+j}|+|c_{n+2j}|\right).
\]
Hence by (8), we get
\begin{eqnarray*}
N^2|c_{n+N}|&\leq&KN\sum\limits_{j=[N/2]}^N\left( |c_{n+j}|+|c_{n+2j}|\right)
\\
&\leq&KN\sum\limits_{j=[N/2]}^{2N}|c_{n+j}|\leq
K\sum\limits_{j=1}^{2N}j|c_{n+j}| \\
&\leq&K\sum\limits_{j=1}^{2N}j\mbox{Re}c_{n+j}\leq 2KNE_n(f).
\end{eqnarray*}
Therefore, we already have
\begin{eqnarray*}
N|c_{n+N}|\leq KE_n(f),\;\;\;N\geq 1,
\end{eqnarray*}
or in other words,
\begin{eqnarray*}
\max\limits_{k\geq 1}k|c_{n+k}|=O(E_n(f)).\hspace{.4in}(9)
\end{eqnarray*}
Analogue to (8), we can also achieve that
\begin{eqnarray*}
2NE_n(f)\geq \sum\limits_{k=1}^Nk\left(\mbox{Re}c_{n+k} +\mbox{Re}%
c_{-n-k}\right).
\end{eqnarray*}
Then the condition $\{c_{n}+c_{-n}\}\in NBVS$, in a similar argument as (9),
leads to
\[
\max\limits_{k\geq 1}k|c_{n+k}+c_{-n-k}|=O(E_n(f)).
\]
Thus
\[
\max\limits_{k\geq 1}|kc_{-n-k}|\leq \max\limits_{k\geq 1}k|c_{n+k}+c_{n-k}|
+\max\limits_{k\geq 1}k|c_{n+k}|=O(E_n(f)).
\]
Lemma 3 is completed. }}

{\normalsize {\rm {\bf Lemma 4 ([8, Lemma 4])}\quad {\it Let $f\in
C_{2\pi},\;c_{n}+c_{-n}\in K(\theta_1),\;n=1,2,\cdots,$ for some $0\leq
\theta_1<\pi/2$. Then
\[
\sum\limits_{k=2n+1}^\infty|c_{k}+c_{-k}|=O(E_n(f)).
\]
} }}

{\normalsize {\rm {\bf Lemma 5.}\quad {\it Let $\left\{c_{n}\right\}_{n=0}^%
\infty\in NBVS$ and $\left\{c_{n}+c_{-n}\right\}_{n=0}^\infty\in NBVS$, then
\[
\left|\sum\limits_{k=1}^n c_{\pm(n+k)}\sin
kx\right|=O\left(\max\limits_{1\leq k\leq
n}k\left(|c_{n+k}|+|c_{-n-k}|\right)\right)
\]
holds uniformly for any $x\in [0,\pi].$} }}

{\normalsize {\rm {\bf Proof.} The case $x=0$ and $x=\pi$ are trivial. When $%
0<x\leq \pi/n$, by the inequality $|\sin x|\leq x$, we have
\begin{eqnarray*}
\left|\sum\limits_{k=1}^n c_{n+k}\sin kx\right|\leq \frac{\pi}{n}%
\sum\limits_{k=1}^nk|c_{n+k}|\leq\pi\max\limits_{1\leq k\leq n}k|c_{n+k}|.%
\hspace{.4in}(10)
\end{eqnarray*}
When $x>\pi/n$, we can find an natural number $m<n$ such that $m\leq
\pi/x<m+1$, then
\begin{eqnarray*}
\left|\sum\limits_{k=1}^m c_{n+k}\sin kx\right|\leq \frac{\pi}{m}%
\sum\limits_{k=1}^m k|c_{n+k}|+\left|\sum\limits_{k=m+1}^n c_{n+k}\sin
kx\right|=:T_1+T_2.
\end{eqnarray*}
It is clear that
\begin{eqnarray*}
|T_1|\leq \pi \max\limits_{1\leq k\leq n}k|c_{n+k}|.\hspace{.4in}(11)
\end{eqnarray*}
By Abel's transformation, we get
\begin{eqnarray*}
T_2=\sum\limits_{k=m+1}^{n-1}\Delta c_{n+k}\sum\limits_{j=1}^k\sin
jx+c_{2n}\sum\limits_{j=1}^n\sin jx-c_{n+m+1}\sum\limits_{j=1}^m\sin jx,
\end{eqnarray*}
hence\footnote{%
Note that $\sum\limits_{j=1}^k\sin jx=O(x^{-1})=O(m+1)$ for $m\leq \pi/x<m+1.
$}
\[
|T_2|\leq C(m+1)\sum\limits_{k=m+1}^{n-1} |\Delta
c_{n+k}|+C\max\limits_{1\leq k\leq n}k|c_{n+k}|.
\]
Note that $\left\{c_{n}\right\}_{n=0}^\infty\in NBVS$. If $n$ is an odd
number, then
\[
\sum\limits_{k=[n/2]}^{n-1}|\Delta
c_{n+k}|=\sum\limits_{k=(n-1)/2}^{n-1}|\Delta c_{n+k}| \leq
K\left(|c_{2n-1}|+|c_{(3n-1)/2}|\right);
\]
if $n$ is an even number, then
\[
\sum\limits_{k=[n/2]}^{n-1}|\Delta
c_{n+k}|\leq\sum\limits_{k=n/2-1}^{n-2}|\Delta c_{n+k}|+ |\Delta c_{2n-1}|
\]
\[
\leq K\left(|c_{2n-2}|+|c_{3n/2-1}|+|c_{2n-1}|+|c_{2n}|\right).
\]
Then, by setting
\[
n^*=:\left\{%
\begin{array}{ll}
\frac{3n-1}{2}, & n\;\mbox{is  odd,} \\
\frac{3n}{2}, & n\;\mbox{is even},%
\end{array}%
\right.
\]
we have
\begin{eqnarray*}
\lefteqn{\su_{k=m+1}^{n-1}|\Delta c_{n+k}|\leq\su_{k=[\log (m+1)/\log
2]+1}^{[\log n/\log 2]}\su_{j=2^{k-1}}^{2^{k}}|\Delta
c_{n+j}|+\su_{k=[n/2]}^{n-1}|\Delta c_{n+k}|}  \nonumber \\
&\leq&K\sum\limits_{k=[\log (m+1)/\log 2]+1}^{[\log n/\log
2]}\left(|c_{n+2^k}|+|c_{n+2^{k+1}}|\right)  \nonumber \\
&&+K\left(|c_{2n-1}|+|c_{n^*}|+ |c_{2n-2}|+|c_{2n-1}|+|c_{2n}|\right)
\nonumber \\
&\leq&K\max\limits_{1\leq k\leq n}k|c_{n+k}|\sum\limits_{k=[\log (m+1)/\log
2]+1}^{[\log n/\log 2]}2^{-k}  \nonumber \\
&&+K\left(|c_{2n-1}|+|c_{n^*}|+ |c_{2n-2}|+|c_{2n-1}|+|c_{2n}|\right)
\nonumber \\
&\leq&K(m+1)^{-1}\max\limits_{1\leq k\leq n}k|c_{n+k}|  \nonumber \\
&&+K\left(|c_{2n-1}|+|c_{n^*}|+ |c_{2n-2}|+|c_{2n-1}|+|c_{2n}|\right).
\end{eqnarray*}
Therefore, with all the above estimates, we deduce that
\begin{eqnarray*}
|T_2|\leq K\max\limits_{1\leq k\leq n}k|c_{n+k}|.\hspace{.4in}(12)
\end{eqnarray*}
Combining (12) with (10) and (11), we finally have
\[
\left|\sum\limits_{k=1}^n c_{n+k}\sin kx\right|= O\left(\max\limits_{1\leq
k\leq n}k|c_{n+k}|\right).
\]
Now write
\[
\sum\limits_{k=1}^n c_{-n-k}\sin
kx=\sum\limits_{k=1}^n\left(c_{n+k}+c_{-n-k}\right)\sin kx
-\sum\limits_{k=1}^n c_{n+k}\sin kx.
\]
Applying the above known estimates, by noting that both $\left\{c_{n}\right%
\}_{n=0}^\infty\in NBVS$ and $\left\{c_{n}+c_{-n}\right\}_{n=0}^\infty\in
NBVS$, we get
\begin{eqnarray*}
\left|\sum\limits_{k=1}^n c_{-n-k}\sin kx\right|&=&O\left(\max\limits_{1\leq
k\leq n}k|c_{n+k}+c_{-n-k}|+\max\limits_{1\leq k\leq n}k|c_{n+k}|\right) \\
&=&O\left(\max\limits_{1\leq k\leq
n}k\left(|c_{n+k}|+|c_{-n-k}|\right)\right).
\end{eqnarray*}
Lemma 5 is proved. }}

{\normalsize {\rm {\bf Lemma 6.}\quad {\it Let $f\in
C_{2\pi},\;\left\{c_{n}\right\}_{n=0}^\infty\in NBVS$ and $%
\left\{c_{n}+c_{-n}\right\}_{n=0}^\infty\in NBVS$, then
\[
\left|\sum\limits_{k=m}^\infty c_{\pm k}\sin
kx\right|=O\left(\max\limits_{k\geq m}k\left( |c_{k}|+|c_{-k}|\right)\right)
\]
holds uniformly for any $x\in [0,\pi]$ and any $m\geq 1$. } }}

{\normalsize {\rm {\bf Proof.} Write
\[
J(x)=\sum\limits_{k=m}^\infty c_{k}\sin kx.
\]
Since for $x=0$ and $x=\pi$, $J(x)=0$, then we may restrict $x$ within $%
(0,\pi)$ without loss of generality. Take $N=[1/x]$ and set\footnote{%
When $N\leq m$, the same argument as in estimating $K_2$ can be applied to
deal with $J(x)$ directly.}
\[
I(x)=\sum\limits_{k=m}^{N-1}c_{k}\sin kx+\sum\limits_{k=N}^\infty c_{k}\sin
kx=:K_1(x)+K_2(x).
\]
Write $\varepsilon_m=\max\limits_{k\geq m}k|c_{k}|.$ It follows from $N=[1/x]
$ that
\[
K_1(x)\leq x\sum\limits_{k=m}^{N-1}k|c_{k}|<x(N-1)\varepsilon_m\leq
\varepsilon_m.
\]
By Abel's transformation, similar to the proof of Lemma 5,
\begin{eqnarray*}
|K_2(x)|&\leq&\left|\sum\limits_{k=N}^\infty\Delta
c_{k}\sum\limits_{\nu=1}^k\sin \nu x-c_{N}\sum\limits_{\nu=1}^{N-1}\sin \nu
x\right| \\
&\leq&\sum\limits_{k=N}^\infty|\Delta c_{k}|\left|\sum\limits_{\nu=1}^k\sin
\nu x\right|+|c_{N}|\left|\sum\limits_{\nu=1}^{N-1}\sin \nu x\right| \\
&\leq&Kx^{-1}\sum\limits_{j=0}^\infty\sum\limits_{k=2^jN}^{2^{j+1}N-1}|%
\Delta c_{k}| +Kx^{-1}|c_{N}| \\
&\leq&K\varepsilon_m\sum\limits_{j=0}^\infty2^{-j}\leq K\varepsilon_m.
\end{eqnarray*}
The treatment of $\sum_{k=m}^{\infty}c_{-k}\sin kx$ is similar. }}

{\normalsize {\rm {\bf Proof of Theorem 4.\quad Necessity.} suppose $f\in
C_{2\pi}.$ From Lemma 3, (1$^{\prime }$) clearly holds, while by Lemma 4, we
see that
\[
\sum\limits_{k=2n+1}^\infty|c_{k}+c_{-k}|=O(E_n(f)),
\]
thus (2$^{\prime }$) holds. }}

{\normalsize {\rm {\bf Sufficiency.} It can be deduced from Theorem 2. }}

{\normalsize {\rm Now, assume that $f\in C_{2\pi }$, by Lemma 3 and Lemma 4,
we see that
\[
\max\limits_{k\geq 1}k|c_{\pm (n+k)}|+\sum\limits_{k=2n+1}^{\infty
}|c_{k}+c_{-k}|=O(E_{n}(f)).
\]%
On the other hand, rewrite $f(x)$ as
\begin{eqnarray*}
f(x) &=&\sum\limits_{k=-2n}^{2n}c_{k}e^{ikx}+i\sum\limits_{k=2n+1}^{\infty
}\left( c_{k}-c_{-k}\right) \sin kx \\
&&+\frac{1}{2}\sum\limits_{k=2n+1}^{\infty }\left( c_{k}+c_{-k}\right)
\left( e^{ikx}+e^{-ikx}\right) ,
\end{eqnarray*}%
then
\[
E_{n}(f)\leq \left\Vert \sum\limits_{k=1}^{n}\left(
c_{n+k}e^{i(n+k)x}+c_{-n-k}e^{-i(n+k)x}\right) \right.
\]%
\[
\left. -\sum\limits_{k=1}^{n}\left(
c_{n+k}e^{i(n-k)x}+c_{-n-k}e^{-i(-n+k)x}\right) \right\Vert
\]%
\[
+\left\Vert \sum\limits_{k=2n+1}^{\infty }\left( c_{k}-c_{-k}\right) \sin
kx\right\Vert +\sum\limits_{k=2n+1}^{\infty }|c_{k}+c_{-k}|
\]%
\[
\leq \left\Vert \sum\limits_{k=1}^{n}c_{n+k}\sin kx\right\Vert +\left\Vert
\sum\limits_{k=1}^{n}c_{-n-k}\sin (-kx)\right\Vert \hspace{2.5cm}
\]%
\[
+\left\Vert \sum\limits_{k=2n+1}^{\infty }\left( c_{k}-c_{-k}\right) \sin
kx\right\Vert +\sum\limits_{k=2n+1}^{\infty }|c_{k}+c_{-k}|.
\]%
Applying Lemma 5 yields that
\[
\left\Vert \sum\limits_{k=1}^{n}c_{n+k}\sin kx\right\Vert +\left\Vert
\sum\limits_{k=1}^{n}c_{-n-k}\sin (-kx)\right\Vert \leq K\max\limits_{1\leq
k\leq n}k\left( |c_{n+k}|+|c_{-n-k}|\right) .
\]%
By noting that $c_{k}-c_{-k}=2c_{k}-(c_{k}+c_{-k})$ and both $\left\{
c_{n}\right\} _{n=0}^{\infty }\in NBVS$ and $\left\{ c_{n}+c_{-n}\right\}
_{n=0}^{\infty }\in NBVS$, we have
\[
\left\Vert \sum\limits_{k=2n+1}^{\infty }\left( c_{k}-c_{-k}\right) \sin
kx\right\Vert =O\left( 2\max\limits_{k\geq 2n+1}k|c_{k}|+\max\limits_{k\geq
2n+1}k|c_{k}+c_{-k}|\right)
\]%
by Lemma 6. Suppose that $\max\limits_{k\geq 2n+1}k|c_{k}|=k_{0}|c_{k_{0}}|.$
Assume that $2n+1\leq k_{0}\leq 4n$, then by the definition of $NBVS$, we
get
\[
|c_{k_{0}}|\leq \sum\limits_{k=k_{0}}^{4n-1}|\Delta c_{k}|+|c_{4n}|\leq
\sum\limits_{k=2n}^{4n}|\Delta c_{k}|+|c_{4n}|\leq K(|c_{2n}|+|c_{4n}|),
\]%
and hence
\[
k_{0}|c_{k_{0}}|\leq K\left( \max\limits_{1\leq k\leq
n}k|c_{n+k}|+4n|c_{4n}|\right) .\hspace{0.4in}(13)
\]%
Meanwhile, for $k_{0}\geq 4n+1$,
\[
k_{0}|c_{k_{0}}|\leq \frac{1}{2}k_{0}\left\vert
c_{k_{0}}-c_{-k_{0}}\right\vert +\frac{1}{2}k_{0}\left\vert
c_{k_{0}}+c_{-k_{0}}\right\vert
\]%
\[
\leq \max\limits_{k\geq 2n+1}k\left\vert c_{k}-c_{-k}\right\vert +O\left(
\sum\limits_{k=2n+1}^{\infty }|c_{k}+c_{-k}|\right) ,
\]%
where the last inequality follows from $k_{0}\geq 4n+1$ and the following
calculation:
\begin{eqnarray*}
|c_{k_{0}}+c_{-k_{0}}| &\leq &\sum\limits_{k=j}^{k_{0}-1}\left\vert \Delta
\left( c_{k}+c_{-k}\right) \right\vert +\left\vert c_{j}+c_{-j}\right\vert
\\
&\leq &\sum\limits_{k=j}^{2j}\left\vert \Delta \left( c_{k}+c_{-k}\right)
\right\vert +\left\vert c_{j}+c_{-j}\right\vert  \\
&\leq &K\left( \left\vert c_{j}+c_{-j}\right\vert +\left\vert
c_{2j}+c_{-2j}\right\vert \right)
\end{eqnarray*}%
for $k_{0}/2+1<j\leq k_{0}$, so that
\[
k_{0}|c_{k_{0}}+c_{-k_{0}}|\leq K\sum\limits_{j=k_{0}/2+1}^{k_{0}}\left(
\left\vert c_{j}+c_{-j}\right\vert +\left\vert c_{2j}+c_{-2j}\right\vert
\right)
\]%
\[
\leq K\sum\limits_{k=2n+1}^{\infty }\left\vert c_{k}+c_{-k}\right\vert .
\]%
By the same technique, for the factor appearing in (13), we also have
\[
4n|c_{4n}|\leq K\left( \max\limits_{k\geq 2n+1}k|c_{k}-c_{-k}|+O\left(
\sum\limits_{k=2n+1}^{\infty }\left\vert c_{k}+c_{-k}\right\vert \right)
\right) .
\]%
Altogether, we see that
\[
\max\limits_{k\geq 2n+1}k|c_{k}|\leq K\max\limits_{1\leq k\leq
n}k|c_{n+k}|+\max\limits_{k\geq 2n+1}k|c_{k}-c_{-k}|+O\left(
\sum\limits_{k=2n+1}^{\infty }\left\vert c_{k}+c_{-k}\right\vert \right)
\]%
holds in any case. With condition that $\{c_{k}+c_{-k}\}\in NBVS$, by a
similar argument we can easily get
\[
\max\limits_{k\geq 2n+1}k\left\vert c_{k}+c_{-k}\right\vert \leq
K\max\limits_{1\leq k\leq n}k\left( |c_{n+k}|+|c_{-n-k}|\right)
\]%
\[
+O\left( \sum\limits_{k=2n+1}^{\infty }\left\vert c_{k}+c_{-k}\right\vert
\right) .
\]%
Combining all the above estimates, we have completed the proof of Theorem 4.
}}

\begin{flushleft}
{\normalsize {\rm {\bf \S 3.4. Proof of Theorem 5 and Theorem 6} }}
\end{flushleft}

{\normalsize {\rm From the definition of $NBVS$, it can be easily deduced
that }}

{\normalsize {\rm {\bf Lemma 8.}\quad {\it Let $\{c_n\}\in NBVS,$ then
\begin{eqnarray*}
\sum\limits_{k=n}^{2n}|\Delta c_k|\log k=O\left(\max\limits_{n\leq k\leq
2n}|c_k|\log k\right),\;n=1,2,\cdots.
\end{eqnarray*}%
} }}

{\normalsize {\rm {\bf Lemma 9 (Xie and Zhou [10]).}\quad {\it Write
\begin{eqnarray*}
\phi_{\pm n}(x):=\sum\limits_{k=1}^n\frac{1}{k}\left(e^{i(k\mp
n)x-e^{-i(k\pm n)x}}\right),
\end{eqnarray*}
then
\[
|\phi_{\pm n}(x)|\leq 6\sqrt{\pi}.
\]%
} }}

{\normalsize {\rm {\bf Proof of Theorem 5.} Denote
\[
\tau_{2n,n}(f,x):=\frac{1}{n}\sum\limits_{k=n}^{2n-1}S_{k}(f,x),
\]
then obviously,
\begin{eqnarray*}
\lim\limits_{n\rightarrow\infty}\|f-\tau_{2n,n}(f)\|_L=0.\hspace{.4in}(14)
\end{eqnarray*}
Write
\[
D_k(x):=\frac{\sin((2k+1)x/2)}{2\sin (x/2)},
\]
\[
D_k^*(x):=\left\{%
\begin{array}{ll}
\frac{\cos(x/2)-\cos((2k+1)x/2)}{2\sin (x/2)} & |x|\leq 1/n, \\
-\frac{\cos((2k+1)x/2)}{2\sin (x/2)} & 1/n\leq |x|\leq\pi,%
\end{array}%
\right.
\]
\[
E_k(x):=D_k(x)+iD_k^*(x).
\]
For $k=n,n+1,\cdots,2n$, we have (see [10], for example)
\begin{eqnarray*}
E_k(\pm x)-E_{k-1}(\pm x)=e^{\pm ikx},\hspace{.4in}(15)
\end{eqnarray*}
\begin{eqnarray*}
E_k(x)+E_k(-x)=2D_k(x), \hspace{.4in}(16)
\end{eqnarray*}
\begin{eqnarray*}
\|E_k\|_L+\|D_k\|_L=O(\log k). \hspace{.4in}(17)
\end{eqnarray*}
By (15) and (16) and applying Abel's transformation, we get
\begin{eqnarray*}
\tau_{2n,n}(f,x)-S_n(f,x)&=&\frac{1}{n}\sum\limits_{k=n+1}^{2n}(2n-k)\left(%
\hat{f}(k)e^{ikx}+ \hat{f}(-k)e^{-ikx}\right)  \nonumber \\
&=&\frac{1}{n}\sum\limits_{k=n+1}^{2n}(2n-k)\left(2\Delta\hat{f}%
(k)D_k(x)-(\Delta\hat{f}(k)- \Delta\hat{f}(-k))E_k(-x)\right)  \nonumber \\
&&+\frac{1}{n}\sum\limits_{k=n}^{2n-1}\left(\hat{f}(k+1)E_k(x)-\hat{f}%
(-k-1)E_k(-x)\right)  \nonumber \\
&&-\left(\hat{f}(n)E_n(x)+\hat{f}(-n)E_n(-x)\right).
\end{eqnarray*}
Thus, (17) and Lemma 8 yield that
\begin{eqnarray*}
\|\tau_{2n,n}(f)-S_n(f)\|_L&=&O\left(\sum\limits_{k=n}^{2n}\big(|\Delta\hat{f%
}(k)|+|\Delta\hat{f}(-k)| \big)\log k+\max\limits_{n\leq |k|\leq 2n-1}|\hat{f%
}(k)|\log k\right) \\
&=&O\left(\max\limits_{n\leq |k|\leq 2n}|\hat{f}(k)|\log k\right).
\end{eqnarray*}
Hence,
\[
\|f-S_n(f)\|_L\leq \|f-\tau_{2n,n}(f)\|_L+O\left(\max\limits_{n\leq |k|\leq
2n}|\hat{f}(k)|\log k\right),
\]
then it follows from (14) and
\begin{eqnarray*}
\lim\limits_{n\rightarrow\infty}\hat{f}(n)\log |n|=0\hspace{.4in}(18)
\end{eqnarray*}
that
\[
\limsup\limits_{n\rightarrow\infty}\|f-S_n(f)\|_L\leq \varepsilon,
\]
in other words,
\begin{eqnarray*}
\lim\limits_{n\rightarrow\infty}\|f-S_n(f)\|_L=0. \hspace{.4in}(19)
\end{eqnarray*}
}}

{\normalsize {\rm Now we come to prove that (19) implies (18). By Lemma 9,
we derive that
\[
\frac{1}{6\sqrt{\pi}}\left|\int_{-\pi}^\pi(f(x)-S_n(f,x))\phi_n(x)dx\right|%
\leq \|f-S_n(f)\|_L,
\]
thus
\begin{eqnarray*}
\sum\limits_{k=1}^n\frac{1}{k}\hat{f}(n+k)=O(\|f-S_n(f)\|_L),
\end{eqnarray*}
especially,
\begin{eqnarray*}
\sum\limits_{k=1}^n\frac{1}{k}\mbox{Re}\hat{f}(n+k)=O(\|f-S_n(f)\|_L).
\end{eqnarray*}
In the same way, we have
\begin{eqnarray*}
\sum\limits_{k=1}^{2n}\frac{1}{k}\mbox{Re}\hat{f}(2n+k)=O(\|f-S_{2n}(f)\|_L).
\end{eqnarray*}
By (3), we deduce that
\begin{eqnarray*}
|\hat{f}(2n)|\log n&\leq& C|\hat{f}(2n)|\sum\limits_{j=1}^{n}\frac{1}{j} \\
&\leq&K({\bf C})\sum\limits_{j=1}^{n}\frac{1}{j}(|\hat{f}(n+j)|+|\hat{f}%
(2n+2j)|) \\
&\leq&K({\bf C})\left(\sum\limits_{j=1}^{n}\frac{1}{j}\mbox{Re}\hat{f}(n+j)+
\sum\limits_{j=1}^{n}\frac{1}{j}\mbox{Re}\hat{f}(2n+2j)\right) \\
&\leq&K({\bf C})\left(\sum\limits_{j=1}^{n}\frac{1}{j}\mbox{Re}\hat{f}(n+j)+
\sum\limits_{j=1}^{2n}\frac{1}{j}\mbox{Re}\hat{f}(2n+j)\right) \\
&\leq&K({\bf C})\Big(\|f-S_n(f)\|_L+\|f-S_{2n}(f)\|_L\Big).
\end{eqnarray*}
Hence, (19) implies
\begin{eqnarray*}
\lim\limits_{n\rightarrow+\infty}|\hat{f}(2n)|\log n=0.
\end{eqnarray*}
The same argument can be applied as well to achieve that
\begin{eqnarray*}
\lim\limits_{n\rightarrow+\infty}|\hat{f}(2n+1)|\log n=0.
\end{eqnarray*}
Therefore, we already have
\[
\lim\limits_{n\rightarrow+\infty}|\hat{f}(n)|\log n=0.
\]
The same discussion leads to
\[
\lim\limits_{n\rightarrow-\infty}|\hat{f}(n)|\log |n|=0.
\]
Theorem 5 is completed. }}

{\normalsize {\rm The proof of Theorem 6 can be proceeded exactly in the
same way as that of Theorem 5 to achieve the required rate. We omit it here.
}}

{\normalsize {\rm \vspace{6mm} }}

\begin{center}
{\normalsize {\rm {\Large {\bf References}} }}
\end{center}

\begin{enumerate}
\item \bf A. S. Belov, {\it On sequential estimate of
best approximation and moduli of continuity by sums of trigonometric
series with quasimonotone coefficients,} \rm  Matem. Zemetik,
51(1992),132-134. (in Russian)

\item \bf T. W. Chaundy and A. E. Jolliffe, {\it The
uniform convergence of a certain class of trigonometric series,} \rm
Proc. London Math. Soc., (15)(1916), 214-216.

\item \bf  R. J. Le and S. P. Zhou, {\it A new condition
for the uniform convergence of certain trigonometric series,} \rm
Acta Math. Hungar., 108(2005),161-169.

\item \bf R. J. Le and S. P. Zhou, {\it  On $L^1$
convergence of Fourier series of complex valued functions,} \rm
Studia Math. Sci. Hungar., to appear.

\item \bf  L. Leindler,{\it On the uniform convergence
and boundedness of a certain class of sine series,} \rm  Anal.
Math., 27(2001), 279-285.

\item \bf L. Leindler, {\it  A new class of numerical
sequences and its applications to sine and cosine series,} \rm
Anal. Math., 28(2002), 279-286.

\item\bf J. R. Nurcombe, {\it  On the uniform convergence
of sine series with quasi-monotone coefficients,} \rm  J. Math.
Anal. Appl., 166(1992), 577-581.

\item\bf W. Xiao, T. F. Xie and S. P. Zhou, {\it The
best approximation rate of certain trigonometric series,} \rm  Ann.
Math. Sinica, 21(A)(2000), 81-88. (in Chinese)

\item \bf  T. F. Xie and S. P. Zhou,{\it On certain
trigonometric series,} \rm Analysis, 14(1994), 227-237.

\item\bf T. F. Xie and S. P. Zhou, {\it $L^1-$%
approximation of Fourier series of complex valued functions,} \rm
Proc. Royal Soc. Edinburg, 126A(1996), 343-353.

\item \bf S. P. Zhou and R. J. Le, {\it  Some remarks on
the best approximation rate of certain trigonometric series,} \rm
Acta Math. Sinica, to appear. (in Chinese)

\item \bf S. P. Zhou and R. J. Le, {\it  A new condition
and applications in Fourier analysis, II,} \rm Advan.
Math.(Beijing), 34(2005), 249-252.

\item\bf A. Zygmund, {\it  Trigonometric Series, 2nd.
Ed., Vol.I,} \rm Cambridge Univ. Press, Cambridge, 1959.
\end{enumerate}

 \vspace{4mm}
\begin{flushleft}
\bf Yu Dansheng\\
\rm Institute of Mathematics\\
 Zhejiang Sci-Tech University\\
 Xiasha Economic Development Area\\
  Hangzhou Zhejiang310018 China\\
  e-mail: danshengyu@yahoo.com.cn

 \vspace{3mm}

\bf Zhou Songping\\
 \rm Institute of Mathematics\\
  Zhejiang
Sci-Tech University \\
Xiasha Economic Development Area\\
 Hangzhou Zhejiang 310018 China\\
 e-mail: szhou@zjip.com
\end{flushleft}

\end{document}